\documentclass[12pt]{iopart}
\usepackage{iopams, graphicx, cite, color}%
\newtheorem{thm}{Theorem}[section]

\newtheorem{lem}[thm]{Lemma}
\newtheorem{rem}[thm]{Remark}

\newcommand{\Real}{\mathbb{R}}
\newcommand{\eps}{\varepsilon}

\newcommand{\rank}{\mathop{\rm rank}}

\newcommand{\abs}[1]{\left\vert#1\right\vert}

\begin{document}
\title[On $\delta'$-like potential scattering]
{On $\delta'$-like potential scattering on star graphs}
\author{S S Man'ko}
\address{Department of Mechanics and Mathematics, Ivan Franko National University of Lviv, 1 Universytetska str., 79000 Lviv, Ukraine }

\ead{s\_\,manko@franko.lviv.ua}

\begin{abstract}
We discuss the potential scattering on the noncompact star graph.
The Schr\"{o}dinger operator with the short-range potential localizing in a neighborhood of the graph vertex is considered.
We study the asymptotic behavior the corresponding scattering matrix in the zero-range limit.
It has been known for a long time that in dimension 1 there is no non-trivial Hamiltonian with the distributional potential $\delta'$, i.e., the $\delta'$ potential acts as a totally reflecting wall. Several authors have, in recent years, studied the scattering properties of the regularizing potentials $\alpha\eps^{-2}Q(x/\eps)$ approximating the first derivative of the Dirac delta function.
A non-zero transmission through the regularized potential has been shown to exist as $\eps\to0$. We extend these results to star graphs
with the point interaction, which is an analogue of $\delta'$ potential on the line.
We prove that generically such a potential on the graph is opaque. We also show that there exists a countable set of resonant intensities for which a partial transmission through the potential occurs.
This set of resonances is referred to as the resonant set and is determined as the spectrum of an auxiliary Sturm-Liouville problem associated with $Q$ on the graph.
\end{abstract}


\section{Introduction}

Schr\"{o}dinger operators with potentials supported on a discrete set of points (such potentials are usually termed ``point interactions'') have attracted considerable attention both in the physical and mathematical literature over several past decades.
One of the reasons for this is that such singular Hamiltonians are widely used in
various application to atomic, nuclear, and solid state physics.
Applications also arise in optics, for instance, in dielectric media where electromagnetic
waves scatter at boundaries or thin layers.
Another reason is that Schr\"{o}dinger operators with point interactions
often form ``solvable'' models in the sense that the resolvents of such operators can
explicitly be calculated. Consequently the spectrum, the eigenfunctions, as well as resonances and scattering quantities, can also be determined explicitly
(see \cite{Albeverio2edition, AlbeverioKurasov} and references therein).
In the physically oriented literature point interactions are often understood as sharply localized potentials, exhibiting a number of interesting features and
unusual effects not seen for regular potentials.

Currently, there is increasing interest in solvable models on graphs in particular, as a reaction to a great deal of progress in fabricating graph-like structures of a semiconductor material, for which graph Hamiltonians represent a natural model (see the survey \cite{Kuchment} for details).
The idea to investigate quantum mechanics of particles confined to a graph
originated with the study of free electron models of organic molecules \cite{Pauling, Platt, RichardBal}.
Among the systems that were successfully modeled by quantum
graphs we mention e.g., single-mode acoustic and
electro-magnetic waveguide networks \cite{FleJohnKu}, Anderson transition \cite{Anderson} and quantum Hall systems with long range
potential \cite{ChaCod}, fracton excitations in fractal structures \cite{AviLuc}, and mesoscopic quantum systems \cite{KSEI}.
The essential component of such graph models is the wavefunction coupling in the vertices.
The interface conditions have to be chosen to make the Hamiltonian
self-adjoint, or in physical terms, to ensure conservation of the probability
current at the vertex.

One of the most natural way to define a graph Hamiltonian corresponding to a point interaction supported by the branching point is to choose an appropriate
self-adjoint extension of the corresponding free Hamiltonian with the interactions points removed.
This approach has been realized in \cite{ExnerSeba:1989}, where the authors found that a vertex joining $n$ graph edges can be described by $n^2$ real parameters defining the interface condition at the vertex.
A general form of such a coupling was described in \cite{KostrykinSch} by a pair of $n\times n$ matrices $\mathrm{A}$, $\mathrm{B}$ such that $\rank{(\mathrm{A},\mathrm{B})}=n$
and $\mathrm{AB}^*$ is self-adjoint;
the boundary values have to satisfy the conditions
\begin{equation}\label{InterfaceCondition}
    \mathrm{A}\Psi(a)+\mathrm{B}\Psi'(a)=0.
\end{equation}
Here the symbol $\Psi(a)$ is used for the column vector of the boundary values at the vertex $a$, and analogously $\Psi'(a)$ stands for the vector of the first derivatives, taken all in the outgoing direction.

Another possible approach to define the vertex coupling is to approximate a quantum graph
with a thin quantum waveguide
(see \cite{AlbeverioCacciapoutiFinco:2007, MoVa, DellTe, ExPo} and references therein).
Different vertex couplings were also discovered in \cite{Exner:1995, Exner:1996, ExnerSeba:1989, ChExTu}.

The Schr\"{o}dinger operators with the Dirac delta-function and its derivatives in potentials have been studying intensively since the eighties of last century (see \cite{Albeverio2edition, AlbeverioKurasov, DemkovOstr} and the references given there). One of the first well-studied Hamiltonian was
the operator with Dirac's delta-function potential $\delta$.
A case of special interest arises when the $\delta$-potential is replaced by the derivative $\delta'$ of the Dirac delta-function.
\v{S}eba~\cite{SebRMP} appears to have been among the first to consider such a Hamiltonian. To define the Hamiltonian
$-\frac{d^2}{d x^2} + \delta'(x)$ he approximated $\delta'$ by regular potentials $\eps^{-2}Q(x/\eps)$ and then investigated the convergence of the corresponding family of regular Schr\"{o}dinger operators
\begin{equation*}
    \mathrm{S}_\eps = -\frac{d^2}{d x^2} + \frac1{\eps^2} Q\Bigl(\frac{x}{\eps}\Bigr).
\end{equation*}
Here~$Q$ is an integrable function.
If $Q$ has zero mean and the integral $\int_\Real t\,Q(t)\,dt$ is equal to $-1$, then the sequence $\eps^{-2}Q(x/\eps)$ converges in the sense of distributions as $\eps\to0$ to $\delta'$, which served as motivation for such an approximation.
\v{S}eba claimed that the operators $\mathrm{S}_\eps$ converge in the
uniform resolvent sense to the direct sum of the unperturbed
half-line Schr\"{o}dinger operators subject to the Dirichlet boundary conditions at $x=0$. From the viewpoint of scattering theory it means that the $\delta'$-barrier is completely opaque, i.e.,  in the limit $\eps\to0$ the potential  $\eps^{-2}Q(x/\eps)$ becomes a totally reflecting wall at $x=0$ splitting the system into two independent subsystems lying on the half-lines $(-\infty,0)$ and $(0,\infty)$.

But this result is in conflict with conclusions reached in \cite{ChristianZolotarIermak03}, where the resonances in the transmission probability for $\delta'$-like potential have been observed.
In that paper an exactly solvable model with a step function $\alpha Q$	was considered. The authors found a discrete set of intensities $\alpha_n$, for which partial transmission
through the limiting $\delta'$-potential occurs.
The values $\alpha_n$ are roots of a transcendent equation depending on
the regularization $Q$.
Exactly solvable models with other piecewise constant potentials as well as nonrectangular
regularizations of $\delta'$  have later been studied in \cite{Zolotaryuk06, ToyamaNogami, Zolotaryuk08, Zolotaryuk10} and the same conclusion has been drawn.

In \cite{GolovatyManko1, GolovatyManko2} the authors approximated the formal Hamiltonian $-\frac{d^2}{d x^2} + q(x)+\alpha\delta'(x)$ on the line by regular Schr\"{o}dinger operators
\begin{equation*}
    \mathrm{S}_\eps(\alpha,Q)=-\frac{d^2}{d x^2} + q(x)+\frac\alpha{\eps^2} Q\Bigl(\frac{x}{\eps}\Bigr)
\end{equation*}
and the similar resonant effect was discovered.
Here the function $Q\in C_0^\infty(-1,1)$ has zero mean, $\alpha$ is a real coupling constant, and $q$ is a real valued fixed potential tending to $+\infty$ as $|x|\to\infty$,
which ensures that the spectrum of $\mathrm{S}_\eps(\alpha,Q)$ is discrete.
The map assigning a selfadjoint operator $\mathrm{S}(\alpha,Q)$ to each pair $(\alpha,Q)$ was constructed there.
The choice of $\mathrm{S}(\alpha,Q)$ was determined by proximity of its energy levels and pure states to those for the Hamiltonian
with regularized potentials for small $\eps$.
It was established that for almost all coupling constants the operator $\mathrm{S}(\alpha,Q)$ is just the direct sum of the Schr\"{o}dinger operators with potential $q$ on half-axes subject to the Dirichlet boundary conditions at the origin.
But in the exceptional case the nontrivial coupling at the origin arises.
The notion of the \textit{resonant set} $\Sigma_Q$, which is the spectrum of the problem
\begin{equation*}
    -g''+\alpha Q(t)g=0,\quad t\in(-1,1),\qquad g'(-1)=g'(1)=0,
\end{equation*}
with respect to the spectral parameter $\alpha$ was introduced in \cite{GolovatyManko1, GolovatyManko2}.
It was also constructed the \textit{coupling function} $\theta_Q:\,\Sigma_Q\to\Real$ defined via $\theta_Q(\alpha)= g_\alpha(1)/g_\alpha(-1)$, where $g_\alpha$ is an eigenfunction corresponding to the eigenvalue $\alpha\in\Sigma_Q$.
In the case when the coupling constant $\alpha$ belongs to the resonance
set, the $\delta'$-barrier admits a partial transmission of a particle, and
an appropriate wave function obeys the interface conditions
\begin{equation*}
    \psi(+0)=\theta_Q(\alpha)\psi(-0),\qquad \theta_Q(\alpha)\psi'(+0)=\psi'(-0).
\end{equation*}

Studies of \cite{GolovatyManko1, GolovatyManko2} have been continued in \cite{GolHryJPA, Manko:2009}. First the findings of \cite{ChristianZolotarIermak03, Zolotaryuk06, ToyamaNogami, Zolotaryuk08, Zolotaryuk10} were generalized in \cite{Manko:2009}, where the scattering on an arbitrary potential of the form $\alpha\eps^{-2}Q(x/\eps)$ has been considered.
It was proved that such a potential is asymptotically transparent only if a coupling constant belongs to the resonant set.
Moreover, the scattering amplitude for the Hamiltonians $\mathrm{S}_\eps(\alpha,Q)$ and $-\frac{d^2}{dx^2}$ converges as $\eps\to0$ to that for the limiting Hamiltonian $\mathrm{S}(\alpha,Q)$.
In \cite{GolHryJPA} the authors have not only pointed out a mistake in \cite{SebRMP}, but also
have showed that the operators $\mathrm{S}_\eps(\alpha,Q)$ converge in the
uniform resolvent sense as $\eps\to0$ to $\mathrm{S}(\alpha,Q)$.
The results of \cite{GolHryJPA, Manko:2009} were obtained for the special case $q=0$; generally, the same can be derived without difficulty.

Many of point interactions on the line were extended to deal with graphs.
One of the best known is the $\delta$ coupling at $n$ edge vertex:
\begin{equation*}
    \psi_1(a)=\ldots=\psi_n(a),\qquad \sum\limits_{l=1}^n\psi'_l(a)
    =\alpha\psi(a)
\end{equation*}
(see, e.g., \cite{Exner:1995, Exner:1996}). Such a model is a generalization of the Hamiltonian
\begin{equation*}
    -\frac{d^2}{dx^2}+\alpha\delta(x),
\end{equation*}
given by $S_{\alpha,\delta}f=-f''$ on the domain
\begin{equation*}
\fl \mathcal{D}(S_{\alpha,\delta})=
    \big\{f\in H^2(\Real\setminus\{0\})\mid f(+0)=f(-0),\quad f'(+0)-f'(-0)=\alpha f(0)
    \}.
\end{equation*}
A similar generalization is possible for the model for Scr\"{o}dinger operators with $\delta'$-interactions, defined by $S_{\beta,\delta'}f=-f''$ on the set of functinctions
\begin{equation*}
\fl \mathcal{D}(S_{\beta,\delta'})=
    \big\{f\in H^2(\Real\setminus\{0\})\mid f'(+0)=f'(-0),\quad f(+0)-f(-0)=\beta f'(0)
    \}.
\end{equation*}
This generalization was described in \cite{Exner:1995, Exner:1996}:
\begin{equation*}
    \psi'_1(a)=\ldots=\psi'_n(a),\qquad \sum\limits_{l=1}^n\psi_l(a)
    =\beta\psi'(a),
\end{equation*}
and was called $\delta'_s$ coupling.
In the following section we briefly sketch the findings of \cite{Manko:2010}, where an analogue for the Hamiltonian $\mathrm{S}_\eps(\alpha,Q)$ was considered on the metric graph.

\subsection{Graph Hamiltonian with the $\delta'$-like potential}
It will be convenient to recall
basic notions of the theory of differential equations on graphs.
By a metric graph $G=(V,E)$ we mean a finite set $V$ of points in $\Real^3$ (whose elements are called \textit{vertices}) together with a set $E$ of smooth regular curves connecting the vertices (the elements of $E$ are called \textit{edges}). The sets of vertices and edges of the graph $G$ are sometimes denoted by $V(G)$ and $E(G)$ respectively.
A map $f:\,G\to\Real$ is said to be a \textit{function on the graph}, and
the restriction of $f$ on the edge $g\in E(G)$ will be denoted by $f_g$.
Let us introduce the space of functions on the graph $C^\infty(G)=\{f \mid
f_g\in C^\infty(\overline{g})\;\hbox{for all}\; g\in E(G)\}$.
Each edge is equipped with a natural parametrization;
the differentiation is performing with respect to the natural parameter.
We denote by $df/dg(a)$ the limit value of the derivative at the point $a\in V(G)$ taken in the direction away from the vertex.
The integral of $f$ over $G$ is the sum of integrals over all edges
\begin{equation*}
    \int_G f\,dG=\sum_{g\in E(G)}\int_g f\,dg.
\end{equation*}
Let $L_2(G)$ be a Hilbert space with the inner product $(u,v)=\int_Gu\overline{v}\,dG$ and with the norm $\|\cdot\|_{L_2(G)}$.
We also introduce the Sobolev spaces $H^k(G)=\{f\in L_2(G) \mid f_g\in H^k(g)\;
\hbox{for all}\; g\in E(G)\}$ and the space of absolutely continuous functions $AC(G)=\{f \mid f_g\in AC(g)\;
\hbox{for all}\; g\in E(G)\}$.

\begin{figure}[hb]
\begin{center}
\includegraphics[scale=1.6]{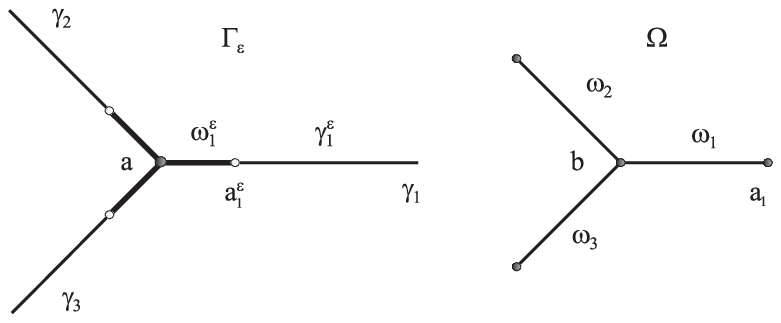}
\end{center}\caption{\label{Graph} Graphs $\Gamma_\eps$ and $\Omega$}
\end{figure}

Let us consider a noncompact star graph $\Gamma$ consisting of three edges $\gamma_1$, $\gamma_2$ and $\gamma_3$. All edges are connected at the vertex $a$.
Then $E(\Gamma)=\{\gamma_1,\gamma_2,\gamma_3\}$ and suppose that all edges are half-lines. We write $a_{j}^\eps$ for the point of intersection of the $\eps$-sphere, centered at $a$, with the edge $\gamma_j\in E(\Gamma)$ and
denote by $\Gamma_\eps$ a sub-partition of $\Gamma$ containing new vertices $a_{1}^\eps$, $a_2^\eps$ and $a_3^\eps$.
Each $a_j^\eps$ divides the edge $\gamma_j$ of the graph $\Gamma$ into two edges $\omega_j^\eps$ and $\gamma_j^\eps$ of the graph $\Gamma_\eps$
(see figure~\ref{Graph}).
Let $\Omega_\eps$ be a star subgraph of $\Gamma_\eps$ such that
$V(\Omega_\eps)=\{a,a_1^\eps,a_2^\eps,a_3^\eps\}$ and
$E(\Omega_\eps)=\{\omega_1^\eps,\omega_2^\eps,\omega_3^\eps\}$.

By $\Omega\subset\Real^3_\xi$ we denote the $\eps^{-1}$-homothety of the graph $\Omega_\eps$, centered at $a$.
Obviously, the resulting graph does not depend on the small parameter $\eps$.
The graph $\Omega$ is a star with a center at the origin $b$ of the auxiliary space $\Real^3_\xi$
and with the vertexes $a_j$,
which are the images of the points $a_j^\eps$.
The edge of $\Omega$ connecting $b$ to $a_j$ will be denoted by $\omega_j$.

We introduce the set $\mathcal{Q}=\big\{Q\in C^\infty(\Omega)\mid \int_\Omega Q\,d\Omega=0\big\}$.
For each nonzero element $Q\in\mathcal{Q}$ let us define the sequence
\begin{equation*}
Q_\eps(x)=
\left\{\begin{array}{ll}
                  \eps^{-2}Q((x-a)/\eps) & \hbox{for $x\in\Omega_\eps$,} \\
                  \phantom{\eps^{-2}V((x-a)/}0 & \hbox{for $x\in\Gamma_\eps\setminus\Omega_\eps$.}
                \end{array}
\right.
\end{equation*}
The potential $Q_\eps$ is referred to as the $\delta'$-like potential.

In \cite{Manko:2010}, it was discovered the family of Schr\"{o}dinger operators on the star graph
\begin{eqnarray*}
\fl   \mathrm{H}_\eps(\alpha,Q)f=-f''+(q +\alpha Q_\eps)f,\\
\fl \mathcal{D}(\mathrm{H}_\eps(\alpha,Q))=\Big\{f\in L_2(\Gamma)\mid
f,\,f'\in AC(\Gamma),\quad -f''+(q+\alpha Q_\eps)f\in L_2(\Gamma),\\
\qquad\qquad\qquad\qquad\qquad f_{\omega_1^\eps}(a)=f_{\omega_2^\eps}(a)=f_{\omega_3^\eps}(a),\quad
\sum\limits_{j=1}^3\frac{d f}{d\omega^\eps_j}(a)=0 \Big\},
\end{eqnarray*}
where $Q\in\mathcal{Q}$, $\alpha$ is a real coupling constant, and
$q$ is a smooth real-valued potential such that
$q_\gamma(x)\to+\infty$ as $\abs{x}\to+\infty$ for all $\gamma\in E(\Gamma)$.
Such behavior of $q$ ensures the discreteness of the spectrum of $\mathrm{H}_\eps(\alpha,Q)$.
A self-adjoint operator $\mathrm{H}(\alpha,Q)$ has been assigned to each pair $(\alpha,Q)$.
The choice of an operator is argued by the close proximity of the energy levels
and the pure states for the Hamiltonians $\mathrm{H}_\eps(\alpha,Q)$ and $\mathrm{H}(\alpha,Q)$ respectively.
Two spectral characteristics of the function $Q$ are introduced in \cite{Manko:2010}: the \textit{resonant set} $\Sigma$, which is the spectrum of the eigenvalue problem
\begin{eqnarray}\label{ResonantSet}
-g''+\alpha Qg=0\quad\hbox{on$\quad\Omega\setminus V(\Omega)$,}
\\\label{ResonantSet1}
\phantom{-}g_{\omega_1}(b)=g_{\omega_2}(b)=g_{\omega_3}(b),\quad\sum\limits_{j=1}^3\frac{d g}{d\omega_j}(b)=0
,
\\\label{ResonantSet2}
\phantom{-}\frac{dg}{d\omega_1}(a_1)=\frac{dg}{d\omega_2}(a_2)=\frac{dg}{d\omega_3}(a_3)
=0,
\end{eqnarray}
and the \textit{coupling function} $\theta:\,\Sigma\to\mathbb{CP}^2$, where $\mathbb{CP}^2$ is the complex projective plane.
The spectrum  of the problem~(\ref{ResonantSet}) is real, discrete and has two accumulation points $\pm\infty$. It consists of simple and double eigenvalues.
Write $\Sigma_1$ for the subset of simple eigenvalues of the problem~(\ref{ResonantSet})
and by $\Sigma_2$ we denote the subset of double eigenvalues.
The resonant set can be represented as the union of $\Sigma_1$ and $\Sigma_2$.

Let $\alpha$ be a simple eigenvalue of the problem~(\ref{ResonantSet}) with an eigenfunction $u_\alpha$ such that $u_\alpha^2(a_1)+u_\alpha^2(a_2)+u_\alpha^2(a_3)=1$, then
introduce
\begin{equation*}
    \theta_1(\alpha)=u_\alpha(a_1),\qquad
    \theta_2(\alpha)=u_\alpha(a_2),\qquad
    \theta_3(\alpha)=u_\alpha(a_3).
\end{equation*}
For $\alpha\in\Sigma_2$ we consider
\begin{eqnarray*} \theta_1(\alpha)=v_\alpha(a_2)w_\alpha(a_3)-v_\alpha(a_3)w_\alpha(a_2)
,\\
\theta_2(\alpha)=v_\alpha(a_3)w_\alpha(a_1)-v_\alpha(a_1)w_\alpha(a_3)
,\\
\theta_3(\alpha)=v_\alpha(a_1)w_\alpha(a_2)-v_\alpha(a_2)w_\alpha(a_1),
\end{eqnarray*}
where $v_\alpha$ and $w_\alpha$ form a base in the corresponding eigenspace such that $\theta_1^2(\alpha)+\theta_2^2(\alpha)+\theta_3^2(\alpha)=1$.
The coupling function is defined via $\theta(\alpha)=(\theta_1(\alpha),\theta_2(\alpha),\theta_3(\alpha))$.
It if easy to check that if we change a base in the eigenspace the point $(\theta_1(\alpha),\theta_2(\alpha),\theta_3(\alpha))\in\mathbb{C}^3$
should be replaced by $(\lambda\theta_1(\alpha),\lambda\theta_2(\alpha),\lambda\theta_3(\alpha))$
with some complex $\lambda$, hence
we shall find it convenient to consider
$\theta$ as a function from the resonant set into $\mathbb{CP}^2$.
In the \textit{non-resonant case}, when the coupling constant $\alpha$ does not belong to the resonant set, $\mathrm{H}(\alpha,Q)$ is the direct sum of the Schr\"{o}dinger operators with the potential $q$ on edges, subject to the Dirichlet boundary conditions at the vertex $a$.
If $\alpha\in\Sigma_1$ (the \textit{case of simple resonance}), then
the operator $\mathrm{H}(\alpha,Q)$ acts via
$\mathrm{H}(\alpha,Q)\psi=-\psi''+q\psi$ on an appropriate set of functions obeying the interface conditions
\begin{equation}\label{InterCond1}
\left.\begin{array}{l}
\theta_2(\alpha)\theta_3(\alpha)\psi_{\gamma_1}(a)=
\theta_1(\alpha)\theta_3(\alpha)\psi_{\gamma_2}(a)=
\theta_1(\alpha)\theta_2(\alpha)\psi_{\gamma_3}(a),\\
\,\qquad\theta_1(\alpha)\frac{d\psi}{d\gamma_1}(a)+
\theta_2(\alpha)\frac{d\psi}{d\gamma_2}(a)+
\theta_3(\alpha)\frac{d\psi}{d\gamma_3}(a)=0
\end{array}
\right.
\end{equation}
that can be expressed by (\ref{InterfaceCondition}) with
\begin{equation*}
\mathrm{A}=
\left(\begin{array}{ccc}
\theta_3(\alpha) & 0 & -\theta_1(\alpha) \\
0 & \theta_3(\alpha) & -\theta_2(\alpha) \\
0 & 0 & 0
\end{array}\right),
\qquad
\mathrm{B}=\left(\begin{array}{ccc}
0 & 0 & 0 \\
0 & 0 & 0 \\
\theta_1(\alpha) & \theta_2(\alpha) & \theta_3(\alpha)
\end{array}\right).
\end{equation*}
In the \textit{case of double resonance} ($\alpha\in\Sigma_2$),
the interface conditions may be written as
\begin{equation}\label{InterCond2}
\left.\begin{array}{l}
\theta_2(\alpha)\theta_3(\alpha)\frac{d\psi}{d\gamma_1}(a)=
\theta_1(\alpha)\theta_3(\alpha)\frac{d\psi}{d\gamma_2}(a)=
\theta_1(\alpha)\theta_2(\alpha)\frac{d\psi}{d\gamma_3}(a),\\
\,\qquad\theta_1(\alpha)\psi_{\gamma_1}(a)+
\theta_2(\alpha)\psi_{\gamma_2}(a)+
\theta_3(\alpha)\psi_{\gamma_3}(a)=0
\end{array}
\right.
\end{equation}
To derive these conditions in the form~(\ref{InterfaceCondition}) we need only interchange the roles of the above matrixes $\mathrm{A}$ and $\mathrm{B}$.

If $\lambda$ is an eigenvalue of the operator $\mathrm{H}(\alpha,Q)$, we will denote by $\mathrm{P}_\lambda$ the orthogonal projector onto the corresponding eigenspace.
Let $\mathrm{P}_\lambda(\eps)$ stand for the orthogonal projector onto the finite dimensional space spanned by all eigenfunctions corresponding to those eigenvalues $\lambda_\eps$ of $\mathrm{H}_\eps(\alpha,Q)$ that $\lambda_\eps\to\lambda$ as $\eps\to0$.
The results of \cite{Manko:2010} may be summed up in the following theorem.
\begin{thm}
All eigenvalues of $\mathrm{H}_\eps(\alpha,Q)$ (except at most a finite number) are bounded as $\eps\to0$.
Let $\lambda_\eps$ be an eigenvalue of $\mathrm{H}_\eps(\alpha,Q)$ bounded as $\eps\to0$,
then $\lambda_\eps$ has a finite limit $\lambda$ that is a point of the spectrum of $\mathrm{H}(\alpha,Q)$. Moreover, $\|\mathrm{P}_\lambda(\eps)-\mathrm{P}_\lambda\|\to0$ as $\eps\to0$.
Conversely, if $\lambda$ is an eigenvalue of $\mathrm{H}(\alpha,Q)$, then there exists an eigenvalue $\lambda_\eps$ of $\mathrm{H}_\eps(\alpha,Q)$ such that $\lambda_\eps\to\lambda$ as $\eps\to0$.
\end{thm}
\smallskip

Although it has been showed the close proximity of the energy levels and pure states for the limiting and regularized Hamiltonians,
we are still in the dark about convergence as $\eps\to0$ of the operators $\mathrm{H}_\eps(\alpha,Q)$ in any topology. As a result, we do not know anything about convergence of the scattering quantities or other physical characteristics.
Our objective in this paper is to give one further motivation for the choice of
the Hamiltonian $\mathrm{H}(\alpha,Q)$.
We shall study the scattering properties of the finite-range potentials $\alpha Q_\eps$ on the graph $\Gamma_\eps$ in the limit $\eps\to0$.
We prove that the scattering coefficients depend on the intensity $\alpha$ and the function $Q$ in such a way that for all values of $\alpha$ the barrier $\alpha Q_\eps$ is completely opaque except for the set $\Sigma$ of resonant values, at which a partial transmission through the potential occurs.
It will also be shown that
the scattering amplitude for the Hamiltonians $\mathrm{H}_\eps(\alpha,Q)$ and $\mathrm{H}_0$ converges as $\eps\to0$ to that for the limiting Hamiltonian $\mathrm{H}(\alpha,Q)$.
Here the Hamiltonian $\mathrm{H}_0$ of a free particle on $\Gamma$  acts via $\mathrm{H}_0f=-f''$ on its domain consisting of those functions from $H^2(\Gamma)$ that are continuous on $\Gamma$ and satisfy the Kirchhoff boundary conditions at the vertex.

\section{Scattering problem for the $\delta'$-like potential on the graph}
Let us state the main result of this paper.
\begin{thm}\label{MainTheorem}
For each $k>0$ and $\alpha\in\mathbb{R}$ the scattering matrix for the operators $\mathrm{H}_\eps(\alpha,Q)$  and $\mathrm{H}_0$ converge as $\eps\to0$ to the scattering matrix for $\mathrm{H}(\alpha,Q)$ and $\mathrm{H}_0$.
\end{thm}
\smallskip

To start with, we briefly treat stationary scattering on the graph $\Gamma$ associated with the Hamiltonians $\mathrm{H}(\alpha,Q)$ and $\mathrm{H}_0$. From now on we assume that $q$ is a zero function, i.e., the operator $\mathrm{H}(\alpha,Q)$ involves no potential and $\mathrm{H}_\eps(\alpha,Q)$ has only perturbed potential $\alpha Q_\eps$.
It is sufficient for us to look at the nontrivial case when the coupling constant $\alpha$ belongs to the resonant set $\Sigma=\Sigma_1\cup\Sigma_2$, since in the opposite case
the particle is always reflected.
Let us introduce a natural parametrization $s\in(0,+\infty)$ on each edge of $\Gamma$ and $\Gamma_\eps$, where $s=0$ corresponds to the vertex $a$.
Consider the incoming monochromatic wave $e^{-iks}$ with $k>0$ coming from infinity along the edge $\gamma_n$. The corresponding wavefunction has the form
\begin{equation}\label{psiN}
    \psi_n(s,k)=
\left\{\begin{array}{ll}
    \phantom{e^{-iks}+}\;\;T_{nm}e^{iks}&\hbox{on $\gamma_m$, $\,m\neq n$,}\\
    e^{-iks}+T_{nn}\,e^{iks}&\hbox{on $\gamma_n$.}
                \end{array}
\right.
\end{equation}
Here $T_{nn}$ are the reflection coefficients, and $T_{nm}$ are the transmission coefficients.
Substituting $\psi_n$ into the matching conditions (\ref{InterCond1}) or (\ref{InterCond2})
we derive the scattering matrix that can be expressed via the coupling function
\begin{equation*}
S(\alpha)=
(-1)^{j-1}2
\left(\begin{array}{ccc}
\theta_1^2(\alpha)-\frac{1}{2}&\theta_1(\alpha)\theta_2(\alpha)&\theta_1(\alpha)\theta_3(\alpha)\\
\theta_1(\alpha)\theta_2(\alpha)&\theta_2^2(\alpha)-\frac{1}{2}&\theta_2(\alpha)\theta_3(\alpha)\\
\theta_1(\alpha)\theta_3(\alpha)&\theta_2(\alpha)\theta_3(\alpha)&\theta_3^2(\alpha)-\frac{1}{2}
\end{array}\right)
\end{equation*}
if $\alpha\in\Sigma_j$ for $j=1,2$.
Note that the scattering matrix does not depend on $k$.

Next let us look in detail at stationary scattering for the Hamiltonians $\mathrm{H}_\eps(\alpha,Q)$ and $\mathrm{H}_0$ and find the limit as $\eps\to0$ of the scattering amplitude.
Consider the incoming monochromatic wave $e^{-iks}$ coming from infinity along the edge $\gamma_3$.

We shall seek a positive-energy solution of the problem
\begin{eqnarray}\label{MainEq}\fl
-y''+\alpha Q_\eps y=k^2 y\quad \hbox{on $\Gamma\setminus \{a\}$},\quad
\phantom{-}y_{\omega_1}(a)=y_{\omega_2}(a)=y_{\omega_3}(a),\quad
\sum_{j=1}^3\frac{dy}{d\omega_j}(a)=0
\end{eqnarray}
that coincides on $\Gamma_\eps\setminus\Omega_\eps$ with $\psi_3$ given by (\ref{psiN}).
Let $u_\alpha=u_\alpha(\xi,\varkappa)$, $v_\alpha=v_\alpha(\xi,\varkappa)$ and $w_\alpha=w_\alpha(\xi,\varkappa)$ be a linearly independent system of solutions of the problem
\begin{equation}\label{ModelEq}
\fl
-g''+\alpha Qg=\varkappa^2g\quad\hbox{on $\Omega\setminus \{b\}$},
\qquad
g_{\omega_1}(b)=g_{\omega_2}(b)=g_{\omega_3}(b),\qquad
\sum_{j=1}^3\frac{dg}{d\omega_j}(b)=0
\end{equation}
on the graph $\Omega$.
Clearly, the functions $u_\alpha((x-a)/\eps,\eps k)$, $v_\alpha=v_\alpha((x-a)/\eps,\eps k)$ and $w_\alpha=w_\alpha((x-a)/\eps,\eps k)$ form a linearly independent system of solutions of the problem (\ref{MainEq}) on $\Omega_\eps$.
Hence the desired solution on $\Omega_\eps$ can be represented as a linear combination of this functions with the coefficients $A$, $B$ and $C$ respectively.
Write $\varkappa=\eps k$.
Since the potential $\alpha Q_\eps$ has discontinuities at the points $a_1^\eps$, $a_2^\eps$ and $a_3^\eps$,  we demand that the conditions
\begin{equation}
\label{MainEq2}
\phantom{-}y_{\gamma_j^\eps}(a_j^\eps)=y_{\omega_j^\eps}(a_j^\eps),\qquad
\frac{dy}{d\gamma_j^\eps}(a_j^\eps)+\frac{dy}{d\omega_j^\eps}(a_j^\eps)=0,\quad
j=1,2,3
\end{equation}
hold.
On substituting our solution into (\ref{MainEq2}) we obtain the linear system
$\mathrm{M}c=r$
for the vector of unknown coefficients $c(\alpha,\varkappa)=(T_{31}, T_{32}, T_{33}, A, B, C)^\top$.
Here $r(\varkappa)=(0,0,0,0,e^{-i\varkappa},-i\varkappa e^{-i\varkappa})^\top$ and
\begin{equation*}\fl \qquad
\mathrm{M}(\alpha,\varkappa)=
\left(\begin{array}{cccccc}
-e^{i\varkappa}&  0&  0&u_\alpha(a_1,\varkappa)&v_\alpha(a_1,\varkappa)&w_\alpha(a_1,\varkappa)\\
-i\varkappa e^{i\varkappa}&  0&  0& \frac{du_\alpha}{d\omega_1}(a_1,\varkappa)&
\frac{dv_\alpha}{d\omega_1}(a_1,\varkappa)&\frac{dw_\alpha}{d\omega_1}(a_1,\varkappa)\\
0& -e^{i\varkappa}&  0& u_\alpha(a_2,\varkappa)&v_\alpha(a_2,\varkappa)&w_\alpha(a_2,\varkappa)\\
0& -i\varkappa e^{i\varkappa}&  0& \frac{du_\alpha}{d\omega_2}(a_2,\varkappa)&
\frac{dv_\alpha}{d\omega_2}(a_2,\varkappa)&\frac{dw_\alpha}{d\omega_2}(a_2,\varkappa)\\
0& 0& -e^{i\varkappa}&  u_\alpha(a_3,\varkappa)&v_\alpha(a_3,\varkappa)&w_\alpha(a_3,\varkappa)\\
0& 0& -i\varkappa e^{i\varkappa}&  \frac{du_\alpha}{d\omega_3}(a_3,\varkappa)&
\frac{dv_\alpha}{d\omega_3}(a_3,\varkappa)&
\frac{dw_\alpha}{d\omega_3}(a_3,\varkappa)
\end{array}\right).
\end{equation*}
Clearly, the functions $u_\alpha(\cdot,0)$, $v_\alpha(\cdot,0)$ and $w_\alpha(\cdot,0)$ form a linearly independent system of solutions of the problem~(\ref{ResonantSet}), (\ref{ResonantSet1}).
In what follows we shall omit the second argument of these functions and keep in mind that it equals zero.
We introduce the determinants
\begin{equation*}
\fl h_0(\alpha)=
\left|\begin{array}{ccc}
u_\alpha(a_1)& v_\alpha(a_1)
& w_\alpha(a_1) \\
u_\alpha(a_2)& v_\alpha(a_2)
& w_\alpha(a_2)\\
u_\alpha(a_3)& v_\alpha(a_3)
& w_\alpha(a_3)
\end{array}\right|,\qquad
h_1(\alpha)=
\left|\begin{array}{ccc}
\frac{du_\alpha}{d\omega_1}(a_1)& \frac{dv_\alpha}{d\omega_1}(a_1)
& \frac{dw_\alpha}{d\omega_1}(a_1) \\
\frac{du_\alpha}{d\omega_2}(a_2)& \frac{dv_\alpha}{d\omega_2}(a_2)
& \frac{dw_\alpha}{d\omega_2}(a_2)\\
\frac{du_\alpha}{d\omega_3}(a_3)& \frac{dv_\alpha}{d\omega_3}(a_3)
& \frac{dw_\alpha}{d\omega_3}(a_3)
\end{array}\right|.
\end{equation*}
Let us also consider the function $h_{nm}$, which is just the determinant $h_{1-n}$ with the $m^{th}$-row  replaced by the $m^{th}$-row of  $h_n$ for $n=0,1$; $m=1,2,3$.

Let us denote by $\Delta(\varkappa,\alpha)$ the determinant of the matrix $\mathrm{M}(\varkappa,\alpha)$.
It admits the asymptotic expansion
\begin{eqnarray*}
\fl
\Delta(\varkappa,\alpha)=h_1(\alpha)+i\varkappa\{3h_1(\alpha)-H_0(\alpha)\}
+\varkappa^2\{-9/2\,h_1(\alpha)+3H_0(\alpha)-H_1(\alpha)\}+O(\varkappa^3)
\end{eqnarray*}
as $\varkappa\to0$, where $H_n(\alpha)=\sum_{m=1}^3h_{nm}(\alpha)$ for $n=0,1$.
Employing Cramer's rule, we derive
\begin{eqnarray}
\label{T1eps}
\fl
\left.\begin{array}{ll}
T_{31}(\varkappa,\alpha)=
\frac{-2i\varkappa(1+i\varkappa)}{\Delta(\varkappa,\alpha)}
&\left|\begin{array}{ccc}
u_\alpha(a_1)& v_\alpha(a_1)& w_\alpha(a_1) \\
\frac{du_\alpha}{d\omega_1}(a_1)& \frac{dv_\alpha}{d\omega_1}(a_1)& \frac{dw_\alpha}{d\omega_1}(a_1) \\
\frac{du_\alpha}{d\omega_2}(a_2)& \frac{dv_\alpha}{d\omega_2}(a_2)& \frac{dw_\alpha}{d\omega_2}(a_2)
\end{array}\right|
\\
&\qquad\qquad-\frac{2\varkappa^2}{\Delta(\varkappa,\alpha)}
\left|\begin{array}{ccc}
u_\alpha(a_1)& v_\alpha(a_1)& w_\alpha(a_1) \\
\frac{du_\alpha}{d\omega_1}(a_1)& \frac{dv_\alpha}{d\omega_1}(a_1)& \frac{dw_\alpha}{d\omega_1}(a_1) \\
u_\alpha(a_2)& v_\alpha(a_2)& w_\alpha(a_2)
\end{array}\right|
+O(\varkappa^3),
\end{array}\right.\\
\nonumber
\fl\left.\begin{array}{ll}
T_{32}(\varkappa,\alpha)=
\frac{2i\varkappa(1+i\varkappa)}{\Delta(\varkappa,\alpha)}
&\left|\begin{array}{ccc}
\frac{du_\alpha}{d\omega_1}(a_1)& \frac{dv_\alpha}{d\omega_1}(a_1)& \frac{dw_\alpha}{d\omega_1}(a_1) \\
u_\alpha(a_2)& v_\alpha(a_2)& w_\alpha(a_2) \\
\frac{du_\alpha}{d\omega_2}(a_2)& \frac{dv_\alpha}{d\omega_2}(a_2)& \frac{dw_\alpha}{d\omega_2}(a_2)
\end{array}\right|
\\
&\qquad\qquad+\frac{2\varkappa^2}{\Delta(\varkappa,\alpha)}
\left|\begin{array}{ccc}
u_\alpha(a_1)& v_\alpha(a_1)& w_\alpha(a_1) \\
u_\alpha(a_2)& v_\alpha(a_2)& w_\alpha(a_2) \\
\frac{du_\alpha}{d\omega_2}(a_2)& \frac{dv_\alpha}{d\omega_2}(a_2)& \frac{dw_\alpha}{d\omega_2}(a_2)
\end{array}\right|
+O(\varkappa^3),
\end{array}\right.
\\
\fl
\nonumber
\left.\begin{array}{ll}
T_{33}(\varkappa,\alpha)=&\big(\Delta(\varkappa,\alpha)\big)^{-1}\Big[
    -h_1(\alpha)+i\varkappa\big\{H_0(\alpha)-2h_{03}(\alpha)-h_1(\alpha)\big\}
    \\
    &-\varkappa^2\big\{H_1(\alpha)+H_0(\alpha)-2(h_{13}(\alpha)+h_{03}(\alpha))
    -1/2h_1(\alpha)\big\}\Big]
+O(\varkappa^3)
\end{array}\right.
\end{eqnarray}
as $\varkappa\to0$.
Similar asymptotics can be derived without difficulty in the case, when the wave package comes from infinity along $\gamma_1$ or $\gamma_2$.
We shall study the asymptotic behavior as $\varkappa\to0$ of the scattering coefficients. Three cases are to be distinguished:
$\Delta(\varkappa,\alpha)=O(1)$,   $\Delta(\varkappa,\alpha)=O(\varkappa)$ and $\Delta(\varkappa,\alpha)=O(\varkappa^2)$ as $\varkappa\to0$.

\begin{lem}\label{LemmaRS}
The set of roots of the equation $h_1(z)=0$ coincides with the resonant set, i.e.,
$h_1$ is a characteristic determinant of the eigenvalue problem
(\ref{ResonantSet})--(\ref{ResonantSet2}).
 \end{lem}
\smallskip

\noindent
\textbf{Proof.}
It is easily seen that the system
$u_\alpha$, $v_\alpha$ and $w_\alpha$ of solutions of the problem~(\ref{ModelEq})
may be chosen so that
\begin{eqnarray*}
u_\alpha(b)=0,\qquad \frac{du_\alpha}{d\omega_1}(b)=0, \qquad
    \frac{du_\alpha}{d\omega_2}(b)=1,\\
v_\alpha(b)=0,\qquad \frac{dv_\alpha}{d\omega_2}(b)=0, \qquad
    \frac{dv_\alpha}{d\omega_1}(b)=1.
\end{eqnarray*}
Clearly, $u_\alpha$ vanishes on $\omega_1$, and  $v_\alpha$ vanishes on $\omega_2$.
Let us consider the linear combination
\begin{equation}\label{combinationG}
g=\frac{dv_\alpha}{d\omega_1}(a_1)\frac{dw_\alpha}{d\omega_2}(a_2)u_\alpha+
\frac{du_\alpha}{d\omega_2}(a_2)\frac{dw_\alpha}{d\omega_1}(a_1)v_\alpha-
\frac{du_\alpha}{d\omega_2}(a_2)\frac{dv_\alpha}{d\omega_1}(a_1)w_\alpha.
\end{equation}
By construction $\frac{dg}{d\omega_1}(a_1)=\frac{dg}{d\omega_2}(a_2)=0$ and $\frac{dg}{d\omega_3}(a_3)=h_1(\alpha)$.
Hence if $h_1(\alpha)=0$ and $g$ is nontrivial on $\Omega$, then $g$  is an eigenfunction of the problem (\ref{ResonantSet}), corresponding to $\alpha$, and therefore $\alpha\in\Sigma$.

Note that at least one of the values $\frac{du_\alpha}{d\omega_j}(a_j)$, $\frac{dv_\alpha}{d\omega_j}(a_j)$ or $\frac{dw_\alpha}{d\omega_j}(a_j)$ is different from zero for $j=1,2,3$.
Now consider the exceptional case, when $g$ is trivial, i.e., all coefficients in (\ref{combinationG}) equal zero. We show then that both values $\frac{du_\alpha}{d\omega_2}(a_2)$ and $\frac{dv_\alpha}{d\omega_1}(a_1)$ are zero.
Conversely, suppose that $\frac{dv_\alpha}{d\omega_1}(a_1)\neq0$. Since $\frac{dv_\alpha}{d\omega_1}(a_1)\frac{dw_\alpha}{d\omega_2}(a_2)=0$, one obtains
$\frac{du_\alpha}{d\omega_2}(a_2)=\frac{dv_\alpha}{d\omega_2}(a_2)=
\frac{dw_\alpha}{d\omega_2}(a_2)=0$, a contradiction. Similarly, $\frac{du_\alpha}{d\omega_2}(a_2)=0$.
It follows that $\frac{du_\alpha}{d\omega_1}(a_1)=\frac{du_\alpha}{d\omega_2}(a_2)=0$ and
$\frac{dv_\alpha}{d\omega_1}(a_1)=\frac{dv_\alpha}{d\omega_2}(a_2)=0$.
Thus the function $g_1=\frac{dv_\alpha}{d\omega_3}(a_3)u_\alpha-\frac{du_\alpha}{d\omega_3}(a_3)v_\alpha$
satisfies $\frac{dg_1}{d\omega_1}(a_1)=\frac{dg_1}{d\omega_2}(a_2)=\frac{dg_1}{d\omega_3}(a_3)=0$.
If $g_1$ is nontrivial on $\Omega$, then it is an eigenfunction of the problem (\ref{ResonantSet}), corresponding to $\alpha$.
If not, then $u_\alpha$ would be a desired eigenfunction.

What is left is to show that $\alpha\in\Sigma$ implies $h_1(\alpha)=0$.
If we choose the new base $u_\alpha$, $v_\alpha$ and $w_\alpha$ such that $u_\alpha$ is an eigenfunction of the problem (\ref{ResonantSet})--(\ref{ResonantSet2})  corresponding to $\alpha$, then the first column of the determinant is zero and
 the assertion follows.
\hfill $\square$
\medskip

\begin{lem}\label{LemmaRS1}
The set of roots of the equation $H_0(z)=0$ that belongs to the resonant set coincides with $\Sigma_2$.
 \end{lem}
\smallskip

\noindent
\textbf{Proof.}
Suppose, contrary to our claim, that
there exists $\alpha\in\Sigma_1$ such that
$H_0(\alpha)=0$.
We choose the new linearly independent system $u_\alpha$, $v_\alpha$ and $w_\alpha$ such that
$u_\alpha$ is an eigenfunction of the problem (\ref{ResonantSet})--(\ref{ResonantSet2}) corresponding to $\alpha$.
Write $V=(\frac{dv_\alpha}{d\omega_1}(a_1),\frac{dv_\alpha}{d\omega_2}(a_2),
\frac{dv_\alpha}{d\omega_3}(a_3))$ and
$W=(\frac{dw_\alpha}{d\omega_1}(a_1),\frac{dw_\alpha}{d\omega_2}(a_2),
\frac{dw_\alpha}{d\omega_3}(a_3))$. These vectors are linearly independent. Indeed, in the opposite case there exists a nonzero number $\lambda$ such that $V=\lambda W$.
Therefore $v_\alpha-\lambda w_\alpha$ is an eigenfunction
of the problem (\ref{ResonantSet})--(\ref{ResonantSet2}) corresponding to $\alpha$ that is linearly independent with $u_\alpha$, a contradiction.

In view of the Lagrange identities
\begin{eqnarray*}
    u_\alpha(a_1)\frac{dv_\alpha}{d\omega_1}(a_1)+
    u_\alpha(a_2)\frac{dv_\alpha}{d\omega_2}(a_2)+
    u_\alpha(a_3)\frac{dv_\alpha}{d\omega_3}(a_3)=0,
\\
    u_\alpha(a_1)\frac{dw_\alpha}{d\omega_1}(a_1)+
    u_\alpha(a_2)\frac{dw_\alpha}{d\omega_2}(a_2)+
    u_\alpha(a_3)\frac{dw_\alpha}{d\omega_3}(a_3)=0
\end{eqnarray*}
the vector $U=(u_\alpha(a_1),u_\alpha(a_2),u_\alpha(a_3))$
is orthogonal to both vectors $V$ and $W$.
The identity $H_0(\alpha)=0$ can be written as
\begin{eqnarray*}
\fl
\left.\begin{array}{l}
u_\alpha(a_1)
\left|\begin{array}{cc}
    \frac{dv_\alpha}{d\omega_2}(a_2)&\frac{dv_\alpha}{d\omega_3}(a_3)\\
    \frac{dw_\alpha}{d\omega_2}(a_2)&\frac{dw_\alpha}{d\omega_3}(a_3)
\end{array}\right|
+
u_\alpha(a_2)
\left|\begin{array}{cc}
    \frac{dv_\alpha}{d\omega_3}(a_3)&\frac{dv_\alpha}{d\omega_1}(a_1)\\
    \frac{dw_\alpha}{d\omega_3}(a_3)&\frac{dw_\alpha}{d\omega_1}(a_1)
\end{array}\right|
\\
\qquad\qquad\qquad\qquad\qquad\qquad\qquad\qquad\qquad
+u_\alpha(a_3)
\left|\begin{array}{cc}
    \frac{dv_\alpha}{d\omega_1}(a_1)&\frac{dv_\alpha}{d\omega_2}(a_2)\\
    \frac{dw_\alpha}{d\omega_1}(a_1)&\frac{dw_\alpha}{d\omega_2}(a_2)
\end{array}\right|=0.
\end{array}\right.
\end{eqnarray*}
It follows that $U$ is orthogonal to the vector product $V\times W$, hence that $U$ is a zero vector, i.e., $u_\alpha(a_1)=u_\alpha(a_2)=u_\alpha(a_3)=0$,  and finally that
$u_\alpha$ is a zero function, which is impossible.

Let $\alpha$ belong to $\Sigma_2$.
We select a linearly independent system $u_\alpha$, $v_\alpha$ and $w_\alpha$ such that
$u_\alpha$ and $v_\alpha$ form a base in the corresponding eigenspace.
Obviously, the determinant $H_0(\alpha)$ vanishes, and the proof is complete.
\hfill $\square$
\medskip

\begin{rem}
The roots of the equation $h_1(z)=0$ as well as $H_0(z)=0$ do not depend on the choice of the linearly independent system of solutions of the problem~(\ref{ModelEq}).
\end{rem}

\begin{lem}\label{LemmaRS2}
The function $H_1$ is different from zero on $\Sigma_2$.
 \end{lem}
\smallskip

\noindent
\textbf{Proof.}
On the contrary, suppose that there exists $\alpha\in\Sigma_2$ such that
$H_1(\alpha)=0$.
Let us select the system
$u_\alpha$, $v_\alpha$ and $w_\alpha$ of linearly independent solutions of the problem~(\ref{ModelEq}) such that
$u_\alpha$ and $v_\alpha$ are eigenfunctions of the problem (\ref{ResonantSet})--(\ref{ResonantSet2}) corresponding to $\alpha$.
Set $U=(u_\alpha(a_1),u_\alpha(a_2),u_\alpha(a_3))$ and
$V=(v_\alpha(a_1),v_\alpha(a_2),v_\alpha(a_3))$. Note that the vectors $U$ and $V$ are linearly independent, since in the opposite case
the functions $u_\alpha$ and $v_\alpha$ would be linearly dependent.

From the Lagrange identities
\begin{eqnarray*}
    u_\alpha(a_1)\frac{dw_\alpha}{d\omega_1}(a_1)+
    u_\alpha(a_2)\frac{dw_\alpha}{d\omega_2}(a_2)+
    u_\alpha(a_3)\frac{dw_\alpha}{d\omega_3}(a_3)=0
\\
    v_\alpha(a_1)\frac{dw_\alpha}{d\omega_1}(a_1)+
    v_\alpha(a_2)\frac{dw_\alpha}{d\omega_2}(a_2)+
    v_\alpha(a_3)\frac{dw_\alpha}{d\omega_3}(a_3)=0
\end{eqnarray*}
it follows that the vector $W=(\frac{dw_\alpha}{d\omega_1}(a_1),\frac{dw_\alpha}{d\omega_2}(a_2),
\frac{dw_\alpha}{d\omega_3}(a_3))$
is orthogonal to $U$ and $V$.
The identity $H_1(\alpha)=0$ can be written in the form
\begin{eqnarray*}
\fl
\frac{dw_\alpha}{d\omega_1}(a_1)
\left|\begin{array}{cc}
    u_\alpha(a_2)&u_\alpha(a_3)\\
    v_\alpha(a_2)&v_\alpha(a_3)
\end{array}\right|
+
\frac{dw_\alpha}{d\omega_2}(a_2)
\left|\begin{array}{cc}
    u_\alpha(a_3)&u_\alpha(a_1)\\
    v_\alpha(a_3)&v_\alpha(a_1)
\end{array}\right|
\\
\qquad\qquad\qquad\qquad\qquad\qquad\qquad
+\frac{dw_\alpha}{d\omega_3}(a_3)
\left|\begin{array}{cc}
    u_\alpha(a_1)&u_\alpha(a_2)\\
    v_\alpha(a_1)&v_\alpha(a_2)
\end{array}\right|=0.
\end{eqnarray*}
Consequently,  $W$ is orthogonal to the vector product $U\times V$.
Analysis similar to that in the proof of the previous lemma gives the assertion of the lemma.
\hfill $\square$
\medskip

We are now in a position to prove the main result.
Without loss of generality we establish the convergence of the transmission coefficient
$T_{31}(\eps k,\alpha)$ as $\eps\to0$.
Other coefficients may be handled in much the same way.

\noindent
\textbf{Proof of Theorem~\ref{MainTheorem}.}
\textit{The non-resonant case.}
Since $\alpha$ is not a resonant coupling constant, $h_1(\alpha)\neq0$.
From (\ref{T1eps}) it immediately follows that $T_{31}(\eps k,\alpha)=O(\eps k)$ as $\eps\to0$.

\textit{The case of simple resonance.}
If $\alpha\in\Sigma_1$, then $h_1(\alpha)=0$ and $H_0(\alpha)\neq0$ in light of Lemmas~\ref{LemmaRS}, \ref{LemmaRS1}. By~(\ref{T1eps})
\begin{equation*}
T_{31}(\eps k,\alpha)=\frac{2}{H_0(\alpha)}
\left|\begin{array}{ccc}
    u_\alpha(a_1)& v_\alpha(a_1)& w_\alpha(a_1) \\
\frac{du_\alpha}{d\omega_1}(a_1)& \frac{dv_\alpha}{d\omega_1}(a_1)& \frac{dw_\alpha}{d\omega_1}(a_1) \\
\frac{du_\alpha}{d\omega_2}(a_2)& \frac{dv_\alpha}{d\omega_2}(a_2)& \frac{dw_\alpha}{d\omega_2}(a_2)
\end{array}\right|
+O(\eps k), \qquad \eps\to0.
\end{equation*}
We continue by choosing the new linearly independent system $u_\alpha$, $v_\alpha$ and $w_\alpha$ such that
$u_\alpha$ is an eigenfunction of the problem (\ref{ResonantSet})--(\ref{ResonantSet2}) corresponding to $\alpha$ satisfying the condition
$u_\alpha^2(a_1)+u_\alpha^2(a_2)+u_\alpha^2(a_3)=1$ and $\frac{dv_\alpha}{d\omega_2}(a_2)=0$.
Observe that
\begin{equation}\label{T10}
T_{31}(0,\alpha)=
2\theta_1(\alpha)\frac{dv_\alpha}{d\omega_1}(a_1)\frac{dw_\alpha}{d\omega_2}(a_2)
H^{-1}_0(\alpha).
\end{equation}
Now let us express $T_{31}(0,\alpha)$ via the coupling function.
Taking into account the Lagrange identities
we see that
\begin{equation*}
\frac{u_\alpha(a_1)}{u_\alpha(a_3)}=
-\frac{\frac{dv_\alpha}{d\omega_3}(a_3)}{\frac{dv_\alpha}{d\omega_1}(a_1)},
\qquad
\frac{u_\alpha(a_2)}{u_\alpha(a_3)}=
\frac{\frac{dv_\alpha}{d\omega_3}(a_3)\frac{dw_\alpha}{d\omega_1}(a_1)-
\frac{dv_\alpha}{d\omega_1}(a_1)\frac{dw_\alpha}{d\omega_3}(a_3)}
{
\frac{dv_\alpha}{d\omega_1}(a_1)\frac{dw_\alpha}{d\omega_2}(a_2)}.
\end{equation*}
Therefore
\begin{eqnarray}\label{h01}\fl
h_{01}(\alpha)=-u_\alpha(a_1)\frac{dv_\alpha}{d\omega_3}(a_3)
\frac{dw_\alpha}{d\omega_2}(a_2)=\frac{u^2_\alpha(a_1)\frac{dv_\alpha}{d\omega_1}(a_1)
\frac{dw_\alpha}{d\omega_2}(a_2)}{\theta_3(\alpha)},\\\label{h02}
\fl
h_{02}(\alpha)=u_\alpha(a_2)\Big\{\frac{dv_\alpha}{d\omega_3}(a_3)
\frac{dw_\alpha}{d\omega_1}(a_1)-\frac{dv_\alpha}{d\omega_1}(a_1)
\frac{dw_\alpha}{d\omega_3}(a_3)\Big\}=\frac{u^2_\alpha(a_2)\frac{dv_\alpha}{d\omega_1}(a_1)
\frac{dw_\alpha}{d\omega_2}(a_2)}{\theta_3(\alpha)}.
\end{eqnarray}
By construction $h_{03}(\alpha)=u^2_\alpha(a_3)\frac{dv_\alpha}{d\omega_1}(a_1)
\frac{dw_\alpha}{d\omega_2}(a_2)/\theta_3(\alpha)$.
Recall that $H_0(\alpha)=h_{01}(\alpha)+h_{02}(\alpha)+h_{03}(\alpha)$.
Combining (\ref{h01}) and (\ref{h02}) yields
\begin{equation*}
    H_0=\frac{\frac{dv_\alpha}{d\omega_1}(a_1)
\frac{dw_\alpha}{d\omega_2}(a_2)}{\theta_3(\alpha)}\{
u^2_\alpha(a_1)+u^2_\alpha(a_2)+u^2_\alpha(a_3)\}=
\frac{\frac{dv_\alpha}{d\omega_1}(a_1)
\frac{dw_\alpha}{d\omega_2}(a_2)}{\theta_3(\alpha)}.
\end{equation*}
Then by
(\ref{T10})
\begin{equation*}
    T_{31}(0,\alpha)=
    2\theta_1(\alpha)\theta_3(\alpha)=T_{31}(\alpha).
\end{equation*}

\textit{The case of double resonance.}
Finally suppose that $\alpha\in\Sigma_2$. By Lemmas~\ref{LemmaRS}, \ref{LemmaRS1} and \ref{LemmaRS2} we have  $h_1(\alpha)=H_0(\alpha)=0$ and $H_1(\alpha)\neq0$, thus
\begin{equation*}
    T_{31}(\eps k,\alpha)=\frac{2}{H_1(\alpha)}
\left|\begin{array}{ccc}
    u_\alpha(a_1)& v_\alpha(a_1)& w_\alpha(a_1) \\
    \frac{du_\alpha}{d\omega_1}(a_1)& \frac{dv_\alpha}{d\omega_1}(a_1)&
    \frac{dw_\alpha}{d\omega_1}(a_1) \\
    u_\alpha(a_2)& v_\alpha(a_2)& w_\alpha(a_2)
\end{array}\right|
+O(\eps k),\qquad \eps\to0
\end{equation*}
by (\ref{T1eps}).
Let us choose the new base $u_\alpha$, $v_\alpha$ and $w_\alpha$ such that $u_\alpha$ and $v_\alpha$ are eigenfunctions of the problem~(\ref{ResonantSet}) corresponding to $\alpha$ satisfying $\theta_1^2(\alpha)+\theta_2^2(\alpha)+\theta_3^2(\alpha)=1$. Consequently,
\begin{equation*}
    T_{31}(0,\alpha)=-2\theta_3(\alpha)\frac{dw_\alpha}{d\omega_1}(a_1)H_1^{-1}(\alpha).
\end{equation*}
Since the Lagrange identities hold
it follows that
\begin{equation*}
    \frac{\frac{dw_\alpha}{d\omega_1}(a_1) }{\theta_1(\alpha)}=
    \frac{\frac{dw_\alpha}{d\omega_2}(a_2) }{\theta_2(\alpha)}=
    \frac{\frac{dw_\alpha}{d\omega_3}(a_3) }{\theta_3(\alpha)}.
\end{equation*}
We thus get
\begin{eqnarray*}\fl\qquad\qquad
h_{11}(\alpha)=\theta_1(\alpha)\frac{dw_\alpha}{d\omega_1}(a_1),\qquad
h_{12}(\alpha)=\theta_2(\alpha)\frac{dw_\alpha}{d\omega_2}(a_2)=
\frac{\theta_2^2(\alpha)}{\theta_1(\alpha)}\frac{dw_\alpha}{d\omega_1}(a_1),\\
\qquad\qquad\qquad h_{13}(\alpha)=\theta_3(\alpha)\frac{dw_\alpha}{d\omega_3}(a_3)=
\frac{\theta_3^2(\alpha)}{\theta_1(\alpha)}\frac{dw_\alpha}{d\omega_1}(a_1).
\end{eqnarray*}
By recalling $H_1(\alpha)=h_{11}(\alpha)+h_{12}(\alpha)+h_{13}(\alpha)$, one obtains
\begin{equation*}H_1(\alpha)=
    \frac{\frac{dw_\alpha}{d\omega_1}(a_1)}{\theta_1(\alpha)}
    \{\theta_1^2(\alpha)+\theta_2^2(\alpha)+\theta_3^2(\alpha)\}=
    \frac{\frac{dw_\alpha}{d\omega_1}(a_1)}{\theta_1(\alpha)},
\end{equation*}
hence
\begin{equation*}
T_{31}(0,\alpha)=-2\theta_1(\alpha)\theta_3(\alpha)
=T_{31}(\alpha),
\end{equation*}
and the proof is complete.
\hfill $\square$
\medskip

\section{An example}
We suppose that the shape of a short-range potential in an actual model can be approximately described as follows
\begin{equation*}
Q=\left\{
                     \begin{array}{ll}
                     \phantom{-}7 & \hbox{ on $\omega_1$, $s\in [0,0.5)$,} \\
                               -7 & \hbox{ on $\omega_1$, $s\in (0.5,1]$,} \\
                     \phantom{-}0 & \hbox{ otherwise.}
                     \end{array}
    \right.
\end{equation*}
We study the scattering properties of the potentials $\alpha\eps^{-2} Q((x-a)/\eps)$ on $\Gamma_\eps$ as $\eps\to0$
for these models. As shown before,
the scattering amplitude for the Hamiltonians $\mathrm{H}_\eps(\alpha,Q)$ and $\mathrm{H}_0$ converges as $\eps\to0$ to that for the limiting Hamiltonian $\mathrm{H}(\alpha,Q)$. Certainly, the resonant set and the coupling function to be found are specific to the given shape $Q$.
\begin{figure}[hb]
\begin{center}
\includegraphics[scale=0.5]{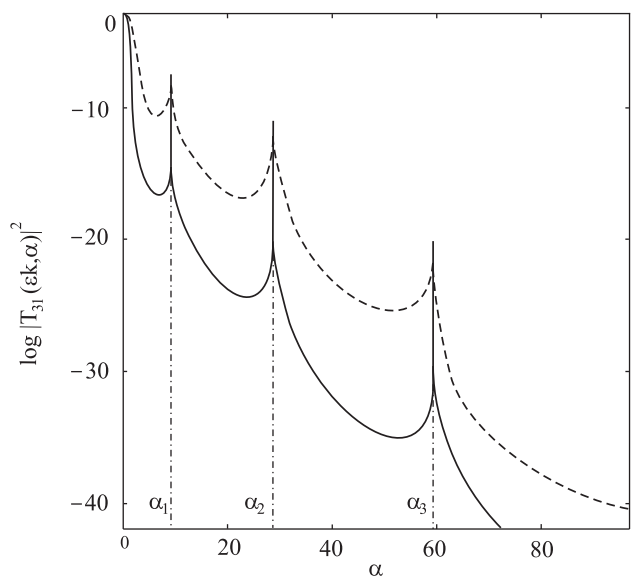}
\end{center}\caption{\label{TransmissionCoeff} Logarithm of the transmission probability $|T_{31}(\eps k,\alpha)|^2$ against $\alpha$ at $\eps k=0.0001$ (solid curve) and 0.01 (dash curve)}
\end{figure}
\begin{table}[h]
\caption{\label{tabone} Resonant intensities and coupling function}
\begin{indented}
\lineup
\item[]\begin{tabular}{*{4}{ccccc}}
\br
$\alpha$&&$\theta_1(\alpha)$&&$\theta_2(\alpha)=\theta_3(\alpha)$\cr
\mr
8.8104 &&$-0.9992$    &&$0.0279$\cr
28.5513&&$0.9999$    &&$0.0012$\cr
59.5701&&$-0.9999$ &&$0.997\times 10^{-4}$\cr
\br
\end{tabular}
\end{indented}
\end{table}
\begin{table}[h]
\caption{\label{tabtwo}Scattering matrix for $\alpha=8.8104$}
\begin{indented}
\lineup
\item[]\begin{tabular}{*{4}{ccccc}}
\br
0.9968 &&$-0.0558$&&$-0.0558$\cr
$-0.0558$&&$-0.9984$&&0.0016\cr
$-0.0558$&& 0.0016&&$-0.9984$ \cr
\br
\end{tabular}
\end{indented}
\end{table}
\begin{table}[h]
\caption{\label{tabthree}Scattering matrix for $\alpha=28.5513$}
\begin{indented}
\lineup
\item[]\begin{tabular}{*{4}{ccccc}}
\br
0.9996 &&$0.0024$&&$0.0024$\cr
$0.0024$&&$-0.9999$&&$0.288\times 10^{-5}$\cr
$0.0024$&& $0.288\times 10^{-5}$&&$-0.9999$ \cr
\br
\end{tabular}
\end{indented}
\end{table}
\begin{table}[h]
\caption{\label{tabfour}Scattering matrix for $\alpha=59.5701$}
\begin{indented}
\lineup
\item[]\begin{tabular}{*{4}{ccccc}}
\br
0.9996 &&$-0.0002$&&$-0.0002$\cr
$-0.0002$&&$-0.9999$&&$0.199\times 10^{-7}$\cr
$-0.0002$&& $0.199\times 10^{-7}$&&$-0.9999$ \cr
\br
\end{tabular}
\end{indented}
\end{table}

Since the function $Q$ vanishes outside the edge $\omega_1$, the resonant set $\Sigma$ coincides with the set $\Sigma_1$ of the simple eigenvalues of (\ref{ResonantSet})--(\ref{ResonantSet2}). In fact, all eigenfunctions are constant outside the edge $\omega_1$ and the resonant set is the spectrum of the problem
\begin{equation*}
-g''+\alpha Qg=0\quad\hbox{on$\quad\omega_1$,}
\qquad
\phantom{-}\frac{d g}{d\omega_1}(a_1)=\frac{d g}{d\omega_1}(b)=0.
\end{equation*}

Table~\ref{tabone} lists the first three positive resonant values of $\alpha$ (numerically computed
using \emph{Maple}) and the corresponding values of the coupling function $\theta(\alpha)$.
Tables~\ref{tabtwo}, \ref{tabthree} and \ref{tabfour} list
the scattering matrixes for different resonant intensities.
We note that the transmission coefficients decay very fast.

The similar effects were observed, when the asymptotic behavior of the scattering coefficients for $\mathrm{H}_\eps(\alpha,Q)$ was discovered directly to the given shape $Q$.
Figure~\ref{TransmissionCoeff} shows the logarithm of the transmission probability $|T_{31}(\eps k,\alpha)|^2$ versus $\alpha$ for $\eps k=0.0001$ and $\eps k=0.01$.
From this figure we see that $|T_{31}(\eps k,\alpha)|^2$ vanishes as $\eps\to0$ except when $\alpha$ takes resonant values which correspond to the spikes of $|T_{31}(\eps k,\alpha)|^2$.
Observe also that the transmission coefficient $T_{31}(\eps k,\alpha)$ decays as $\alpha\to\infty$.

\ack The author expresses his gratitude to Prof Yu~Golovaty for stimulating and helpful discussions.

\end{document}